# ON THE SECOND MOMENT OF THE NUMBER OF CROSSINGS BY A STATIONARY GAUSSIAN PROCESS

By Marie F. Kratz and José R. León[1]

*MAP5 UMR 8145 and SAMOS-MATISSE UMR 8595, and UCV Caracas*

Cramér and Leadbetter introduced in 1967 the sufficient condition

$$\frac{r''(s) - r''(0)}{s} \in L^1([0,\delta], dx), \qquad \delta > 0,$$

to have a finite variance of the number of zeros of a centered stationary Gaussian process with twice differentiable covariance function $r$. This condition is known as the Geman condition, since Geman proved in 1972 that it was also a necessary condition. Up to now no such criterion was known for counts of crossings of a level other than the mean. This paper shows that the Geman condition is still sufficient and necessary to have a finite variance of the number of any fixed level crossings. For the generalization to the number of a curve crossings, a condition on the curve has to be added to the Geman condition.

**1. Introduction and main result.** Let $X = \{X_t,\ t \in \mathbb{R}\}$ be a centered stationary Gaussian process. Its correlation function $r$ is supposed to be twice differentiable and to satisfy on $[0, \delta]$, with $\delta > 0$,

$$r(\tau) = 1 + \frac{r''(0)}{2}\tau^2 + \theta(\tau)$$

(1)

$$\text{with } \theta(\tau) > 0, \frac{\theta(\tau)}{\tau^2} \to 0, \frac{\theta'(\tau)}{\tau} \to 0, \theta''(\tau) \to 0, \text{ as } \tau \to 0.$$

The nonnegative function $L$ defined by $\theta''(\tau) := \tau L(\tau)$ will be referred to as the Geman function.

Let us consider a continuous differentiable real function $\psi$ and let us define, as in [2], the number of crossings of the function $\psi$ by the process

Received December 2004; revised June 2005.
[1]Supported in part by the Project No. 97003647 "Modelaje Estocástico Aplicado" of the Agenda Petróleo of FONACIT Venezuela.

$X$ on an interval $[0,t]$ ($t \in \mathbb{R}$), as the random variable $N_t^\psi = N_t(\psi) = \#\{s \leq t : X_s = \psi_s\}$.

The number $N_t^\psi$ of $\psi$-crossings by $X$ can also be seen as the number of zero crossings $N_t^Y(0)$ by the nonstationary (but stationary in the sense of the covariance) Gaussian process $Y = \{Y_s,\ s \in \mathbb{R}\}$, with $Y_s := X_s - \psi_s$, that is, $N_t^\psi = N_t^Y(0)$.

Regarding the moments of the number of crossings by $X$, one of the most well-known first results was obtained by Rice [9] for a given level $x$, namely

$$\mathbb{E}[N_t(x)] = te^{-x^2/2}\sqrt{-r''(0)}/\pi.$$

This equality was proved two decades later by Itô [7] and Ylvisaker [11], providing a necessary and sufficient condition to have a finite mean number of crossings:

$$\mathbb{E}[N_t(x)] < \infty \iff -r''(0) < \infty.$$

Also in the 1960s, following on the work of Cramér, generalization to curve crossings and higher-order moments for $N_t(\cdot)$ were considered in a series of papers by Cramér and Leadbetter [2] and Ylvisaker [12].

Moreover, Cramér and Leadbetter [2] provided an explicit formula for the second factorial moment of the number of zeros of the process $X$, and proposed a sufficient condition on the correlation function of $X$ in order to have the random variable $N_t(0)$ belonging to $L^2(\Omega)$, namely

$$\text{If } L(t) := \frac{r''(t) - r''(0)}{t} \in L^1([0,\delta], dx) \qquad \text{then } \mathbb{E}[N_t^2(0)] < \infty.$$

Geman [6] proved that this condition was not only sufficient but also necessary:

(2) $\quad \mathbb{E}[N_t^2(0)] < \infty \iff L(t) \in L^1([0,\delta], dx) \qquad$ (Geman condition).

This condition held only when choosing the level as the mean of the process.

Generalizing this result to any given level $x$ and to some differentiable curve $\psi$ has been subject to some investigation and nice papers, such as the ones of Cuzick [4, 5] proposing sufficient conditions. But to get necessary conditions remained an open problem for many years. The solution of this problem is enunciated in the following theorem.

THEOREM.

(1) *For any given level $x$, we have*

$$\mathbb{E}[N_t^2(x)] < \infty \iff \exists \delta > 0, L(t) = \frac{r''(t) - r''(0)}{t} \in L^1([0,\delta], dx)$$

(*Geman condition*).

(2) *Suppose that the continuous differentiable real function $\psi$ is such that*

(3)
$$\exists \delta > 0 \qquad \int_0^\delta \frac{\gamma(s)}{s}\, ds < \infty$$

*where $\gamma(\cdot)$ is the modulus of continuity of $\dot\psi$.*

*Then*

$$\mathbb{E}[N_t^2(\psi)] < \infty \iff L(t) \in L^1([0,\delta], dx).$$

REMARK. This smooth condition on $\psi$ is satisfied by a large class of functions which includes in particular functions whose derivatives are Hölder.

Finally let us mention the work of Belyaev [1] and Cuzick [3, 4, 5] who proposed some sufficient conditions to have the finiteness of the $k$th (factorial) moments for the number of crossings for $k \geq 2$. When $k \geq 3$, the difficult problem of finding necessary conditions when considering levels other than the mean is still open.

**2. Proof.** Generalizing the formula of Cramér and Leadbetter ([2], page 209) concerning the zero crossings, the second factorial moment $M_2^\psi$ of the number of $\psi$-crossings can be expressed as

(4)
$$M_2^\psi = \int_0^t \int_0^t \int_{R^2} |\dot x_1 - \dot\psi_{t_1}||\dot x_2 - \dot\psi_{t_2}|$$
$$\times p_{t_1,t_2}(\psi_{t_1}, \dot x_1, \psi_{t_2}, \dot x_2)\, d\dot x_1\, d\dot x_2\, dt_1\, dt_2,$$

where $p_{t_1,t_2}(x_1, \dot x_1, x_2, \dot x_2)$ is the density of the vector $(X_{t_1}, \dot X_{t_1}, X_{t_2}, \dot X_{t_2})$ that is supposed nonsingular for all $t_1 \neq t_2$. The formula holds whether $M_2^\psi$ is finite or not.

We also have

(5)
$$M_2^\psi = 2\int_0^t \int_{t_1}^t p_{t_1,t_2}(\psi_{t_1}, \psi_{t_2})$$
$$\times \mathbb{E}[|\dot X_{t_1} - \dot\psi_{t_1}||\dot X_{t_2} - \dot\psi_{t_2}|\,|\,X_{t_1} = \psi_{t_1}, X_{t_2} = \psi_{t_2}]\, dt_2\, dt_1,$$

where $p_{t_1,t_2}(x_1, x_2)$ is the density of $(X_{t_1}, X_{t_2})$.

From now on, let us put $t_2 = t_1 + \tau$, $\tau > 0$.

The method used to prove that the Geman condition keeps being the sufficient and necessary condition to have $M_2^\psi$ finite can be sketched into three steps.

The first one consists in using the following regression model to compute the expectation in $M_2^\psi$:

$$(R) \qquad \begin{aligned} \dot{X}_{t_1} &= \zeta + \alpha_1(\tau) X_{t_1} + \alpha_2(\tau) X_{t_1+\tau}, \\ \dot{X}_{t_1+\tau} &= \zeta^* - \beta_1(\tau) X_{t_1} - \beta_2(\tau) X_{t_1+\tau}, \end{aligned}$$

where $(\zeta, \zeta^*)$ is jointly Gaussian such that

$$(6) \qquad \mathrm{Var}(\zeta) = \mathrm{Var}(\zeta^*) := \sigma^2(\tau) = -r''(0) - \frac{r'^2(\tau)}{1 - r^2(\tau)},$$

$$(7) \qquad \rho(\tau) := \frac{\mathrm{Cov}(\zeta, \zeta^*)}{\sigma^2(\tau)} = \frac{-r''(\tau)(1 - r^2(\tau)) - r'^2(\tau) r(\tau)}{-r''(0)(1 - r^2(\tau)) - r'^2(\tau)},$$

and where

$$\alpha_1 = \alpha_1(\tau) = \frac{r'(\tau) r(\tau)}{1 - r^2(\tau)}; \qquad \alpha_2 = \alpha_2(\tau) = -\frac{r'(\tau)}{1 - r^2(\tau)}$$

$$\beta_1 = \beta_1(\tau) = \alpha_2(\tau); \qquad \beta_2 = \beta_2(\tau) = \alpha_1(\tau).$$

In the second step, the expectation, formulated in terms of $\zeta$ and $\zeta^*$, will be expand into Hermite polynomials. Recall that the Hermite polynomials $(H_n)_{n \geq 0}$, defined by $H_n(x) = (-1)^n e^{x^2/2} \frac{d^n}{dx^n} e^{-x^2/2}$, constitute a complete orthogonal system in the Hilbert space $L^2(\mathbb{R}, \varphi(u) \, du)$, $\varphi$ denoting the standard normal density.

Finally, this Hermite expansion will allow us to find, in an easier way, lower and upper bounds for $M_2^\psi$. Nevertheless, it will required a fine study in the neighborhood of 0, on one hand on the correlation function $r$ of $X$ and its derivatives, showing in particular the close relation between the existence of the Geman function $L$ and the existence of $r^{(iv)}(0)$, on the other hand, on the correlation function $\rho$ of the r.v. $\zeta$ and $\zeta^*$ of the model $(R)$. It will be presented in the two first lemmas below. Moreover, since the bounds will be expressed in terms of the variance $\sigma^2(\tau)$ of the r.v. $\zeta$ (or $\zeta^*$), an interesting lemma (see Lemma 3 below) will show that the behavior of $L$ is closely related to the behavior of $\sigma^2(\tau)$.

LEMMA 1.

(i) If $r^{(iv)}(0) = +\infty$, then $\lim_{\tau \to 0} \frac{L(\tau)}{\tau} = +\infty$.
(ii) If $r^{(iv)}(0) < +\infty$, then $\lim_{\tau \to 0} \frac{L(\tau)}{\tau} = \frac{r^{(iv)}(0)}{2}$.

LEMMA 2. *For $\tau$ belonging to a neighborhood of 0:*

(i) $|\frac{r'(\tau)}{\sigma(\tau)}|$ *is bounded;*
(ii) $\rho(\tau) \leq 0$.

LEMMA 3. *For $\tau$ belonging to a neighborhood of 0:*

(i) $\frac{\sigma^2(\tau)}{\tau} \leq L(\tau) \leq (2+C)\frac{\sigma^2(\tau)}{\tau}$, *with $C \geq 0$;*

(ii) *For $\delta > 0$,* $\int_0^\delta \frac{\sigma^2(\tau)}{\sqrt{1-r^2(\tau)}}\, d\tau < \infty \Leftrightarrow \int_0^\delta L(\tau)\, d\tau < \infty$ *(Geman condition).*

The proofs of the lemmas are given in [8].

To illustrate the method, we will present the complete proof when considering a fixed level $x$. For the case of curve-crossings, you can refer to [8].

So suppose $\dot{\psi}_s = 0$ and $\psi_s \equiv x, \forall s$.

Let $C$ be a positive constant which may vary from equation to equation. By using the regression $(R)$, $M_2^x$ can be written as

$$M_2^x = 2\int_0^t (t-\tau) p_\tau(x,x) \sigma^2(\tau) A(m,\rho,\tau)\, d\tau,$$

where

$$A(m,\rho,\tau) := \mathbb{E}\left|\left(\frac{\zeta}{\sigma(\tau)} + \frac{r'(\tau)}{(1+r(\tau))\sigma(\tau)}x\right)\left(\frac{\zeta^*}{\sigma(\tau)} - \frac{r'(\tau)}{(1+r(\tau))\sigma(\tau)}x\right)\right|,$$

and $p_\tau(x,x) := p_{0,\tau}(x,x)$.

Note that

$$(8) \quad M_2^x \geq M_2^{x,\delta} := 2\int_0^\delta (t-\tau) p_\tau(x,x) \sigma^2(\tau) A(m,\rho,\tau)\, d\tau, \qquad \delta \in [0,\tau].$$

Now, by using Mehler's formula (see, e.g., [10]), we have

$$A(m,\rho,\tau) = \sum_{k=0}^\infty a_k(m) a_k(-m) k! \rho^k(\tau) \qquad \text{where } m = m(\tau) := \frac{r'(\tau)x}{(1+r(\tau))\sigma(\tau)},$$

$|m| = |m(\tau)|$ being bounded because of (i) of Lemma 2, and $a_k(m)$ are the Hermite coefficients of the function $|\cdot - m|$, given by

$$a_0(m) = \mathbb{E}|Z-m| \qquad Z \text{ being a standard Gaussian r.v.}$$
$$= m[2\Phi(m)-1] + \sqrt{\frac{2}{\pi}} e^{-m^2/2},$$
$$a_1(m) = (1-2\Phi(m)) = -\sqrt{\frac{2}{\pi}} \int_0^m e^{-u^2/2}\, du$$

and

$$a_l(m) = \sqrt{\frac{2}{\pi}} \frac{1}{l!} H_{l-2}(m) e^{-m^2/2}, \qquad l \geq 2.$$

Let us show that $M_2^x < \infty$ under the Geman condition.

Since by Cauchy–Schwarz inequality

$$|A(m,\rho,\tau)| \leq \sum_{k=0}^{\infty} |a_k(m)a_k(-m)|k! \leq (\mathbb{E}[(Y-m)^2]\mathbb{E}[(Y+m)^2])^{1/2},$$

with $Y$ a standard normal r.v., there follows

$$M_2^x \leq I_2 := 2\int_0^t (t-\tau)p_\tau(x,x)\sigma^2(\tau)(a_0(m)a_0(-m)+1+m^2)\,d\tau.$$

Hence, $m^2$ being bounded, we obtain $I_2 \leq C\int_0^t(t-\tau)p_\tau(x,x)\sigma^2(\tau)\,d\tau$.

The study of this last integral reduces to the one on $[0,\delta]$ because of the uniform continuity outside of a neighborhood of 0, so we can conclude that it is finite if $L \in L^1[0,\delta]$, by using Lemma 3(ii).

Let us look now at the reverse implication.

Suppose that $M_2^x < \infty$, and so, via (8), that $M_2^{x,\delta} < \infty$.

Let us compute $A(m,\rho,\tau)$ and bound it below.

By using the parity of the Hermite polynomials and the sign of $\rho$ given in (ii) of Lemma 2, we obtain

$$A(m,\rho,\tau) = a_0^2(m) + |\rho(\tau)|a_1^2(m) + \sum_{k=1}^{\infty} a_{2k}^2(m)(2k)!\rho^{2k}(\tau)$$

$$+ |\rho|\sum_{k=1}^{\infty} a_{2k+1}^2(m)(2k+1)!\rho^{2k}(\tau)$$

$$\geq a_0^2(m) = \left(-ma_1(m) + \sqrt{\frac{2}{\pi}}e^{-m^2/2}\right)^2$$

$$\geq \frac{2}{\pi}e^{-m^2} \geq C \quad \text{(since } |m| < \infty\text{)}.$$

Hence

$$M_2^{x,\delta} \geq C\int_0^\delta (t-\tau)p_\tau(x,x)\sigma^2(\tau)\,d\tau \geq C\int_0^\delta (t-\tau)\frac{\sigma^2(\tau)}{\sqrt{1-r^2(\tau)}}\,d\tau.$$

An application of Lemma 3(ii), yields that $M_2^{x,\delta} < \infty$ implies the Geman condition.

The proof of the general case follows the same approach. It requires also to use Taylor formula for $\psi$ and to introduce the modulus of continuity of $\dot\psi$ to express the expectation in the integrand of $M_2^\psi$ into two terms, one on which will be applied the described method, the other related to the modulus of continuity of $\dot\psi$, which is bounded thanks to the condition (3) of the theorem (for more details, see [8]).

**Acknowledgments.** J. León is grateful to the SAMOS-MATISSE (Univ. Paris 1) for their invitation in October 2004. We wish to thank the referees for some very useful comments.

U.F.R. DE MATHÉMATIQUES ET INFORMATIQUE
UNIVERSITÉ RENÉ DESCARTES, PARIS V
45 RUE DES SAINTS-PÈRES
75270 PARIS CEDEX 06
FRANCE
E-MAIL: kratz@math-info.univ-paris5.fr
URL: www.math-info.univ-paris5.fr/~kratz

ESCUELA DE MATEMÁTICA
FACULTAD DE CIENCIAS
UNIVERSIDAD CENTRAL DE VENEZUELA
A.P. 47197 LOS CHAGUARAMOS
CARACAS 1041-A
VENEZUELA
E-MAIL: jleon@euler.ciens.ucv.ve
URL: euler.ciens.ucv.ve/